\documentclass{ws-procs9x6}

\usepackage{amsmath}
\usepackage{amssymb}
\usepackage{graphicx}
\def\be{\begin{equation}}
\def\ee{\end{equation}}
\def\bea{\begin{eqnarray}}
\def\eea{\end{eqnarray}}
\def\p{\partial}

\newcommand{\R}{{\mathbb R}}
\newcommand{\pp}{{\mathbb P}}

\begin{document}

\title{OPTIMAL COMBINATION OF DATA MODES IN INVERSE PROBLEMS: MAXIMUM
COMPATIBILITY ESTIMATE
}

\author{Mikko Kaasalainen}
\address{
Department of Mathematics\\
Tampere University of Technology\\
P.O. Box 553, 33101 Tampere\\ 
Finland}

\begin{abstract}
We present an optimal strategy for the relative weighting of different data 
modes in inverse problems, and derive the maximum compatibility 
estimate (MCE) that corresponds to the maximum likelihood or maximum
a posteriori estimates in the case of a single data mode. MCE is not
explicitly dependent on the noise levels, scale factors or numbers of data
points of the complementary data modes, and can be determined without the
mode weight parameters. As a case study, we consider the problem of
reconstructing the shape of a body in $\R^3$ from the boundary curves
(profiles) and volumes (brightness values) of its generalized projections.
\end{abstract}

\keywords{Inverse problems, computational geometry, three-dimensional 
polytopes}

\bodymatter
\section{Introduction}
In many inverse problems, various complementary data modes are available.
For example, constructing the shape model of a body in $\R^3$
is typically based on projectionlike data at various viewing geometries. 
In this paper, we consider the case where images ${\cal I}(\omega,\omega_0)$
(generalized projections)
obtained at viewing and illumination directions 
$\omega,\omega_0\in S^2$ are available, but the reliable infomation in these
images is only contained in the boundary curves $\p\cal I$
between the dark background 
or a shadow and the illuminated portion of the target surface. This is
a typical case in adaptive optics data in astrophysics, where
the coverage of viewing geometries is also seldom wide enough to enable a
full reconstruction of the model from images alone \cite{carry}. Thus we 
include the possibility of augmenting the image dataset with a set of
measured brightnesses (volumes of the generalized projections)
$L(\omega,\omega_0)$ of the target at
various observing geometries. 

\section{Case study: generalized projections}

We consider the inverse problem of determining the shape
of a body ${\cal B}\in\R^3$ 
from some measured profiles of generalized projections 
$\p{\cal I}(\omega_i,\omega_{0i})$, $i=1,\ldots,n$
and their volumes $L(\omega_{0i},\omega_i)$, $i=1,\ldots,m$\cite{ipi,genproj}.

Our goal is to construct a total goodness-of-fit measure $\chi^2_{\rm tot}$ 
\be
\chi^2_{\rm tot}=\chi^2_L+\lambda_\p\chi^2_\p+\lambda_R g(P),
\ee
where $L$ denotes brightness data, $\p$ generalized profiles, and $R$ 
regularizing functions $g(P)$ (see Ref.~\refcite{ipi} for discussion 
of these), where $P\in\R^p$ is the vector of
model parameters. Determining an optimal value for $\lambda_\p$ 
(and $\lambda_R$) is part of the inverse problem.

The volumes of generalized projections are also called total or
disk-integrated brightnesses \cite{genproj}:
\be
L(\omega_0,\omega)=\int_{{\cal A}_+} 
R(x;\omega_0,\omega)\langle\omega,\nu(x)\rangle\,  d\sigma(x),
\ee
where ${\cal A}_+$ is the set of visible and illuminated points $x\in\cal B$
\cite{genproj},
$\nu(x)\in S^2$ and $d\sigma(x)$ are, respectively, the
outward surface normal and surface patch of $\cal B$,
and $R(x;\omega_0,\omega)\in \R$ describes the intensity of 
scattered light at the point $x$ on the surface.
In its basic form,
\be
\chi^2_L=\sum_{i} [L^{\rm (obs)}(\omega_{0i},\omega_i)
-L^{\rm (mod)}(\omega_{0i},\omega_i)]^2\label{lcchisq}
\ee
(assuming a constant noise level; 
see Ref.~\refcite{genproj} and references therein
for modifications and variations of this). $L$-data
on $S^2\times S^2$ uniquely determine a convex body and the
solution is stable \cite{genproj}, but $L$-data do not carry information
on nonconvexities in most realistically available $S^2\times S^2$ 
geometries in practice. 

For many typical adaptive optics targets in our
solar system, the generalized profiles are starlike
due to the proximity of $\omega$ and $\omega_0$ and some regularity
of the target shape at global scale \cite{carry}. Then 
we can write $\chi^2_\p$ by considering, for each profile $i$,
their observed and modelled
maximal radii (from some point $\varkappa_0\in\R^2$ within the profile) 
on the projection plane at direction angles 
$\alpha_{ij}$ (starting from a chosen
coordinate direction):
\be
\chi^2_\p=\sum_{ij} [r_{\rm max}^{\rm (obs)}(\alpha_{ij})
-r_{\rm max}^{\rm (mod)}(\alpha_{ij})]^2.\label{aochisq}
\ee

We now represent the body $\cal B$ as a polytope. Let
two vertices $a$ and $b$ of a facet have projection points $\varkappa_a$,
$\varkappa_b$.
The intersection point $\varkappa$ 
of the radius line at $\alpha$ and the projection of the
facet edge $ab$ is readily determined.
The model $r_{\rm max}(\alpha)$ 
can now be determined by going through all eligible
facet edges and their intersection points $\varkappa_{ab}(\alpha)$:
\be
r_{\rm max}^{\rm (mod)}(\alpha)=\max\Big\{\Vert
\varkappa_{ab}(\alpha)-\varkappa_0\Vert
\Big\vert a,b\in {\cal V}_+\Big\},
\ee
where ${\cal V}_+$ is the set of vertices of the set of facets 
$\tilde{\cal A}_+$ approximating ${\cal A}_+$. The set $\tilde{\cal A}_+$ 
is determined by ray-tracing \cite{genproj}. In general, facet edge circuits 
$\p\tilde{\cal A}_+$ approximating $\p{\cal A}_+$ (and corresponding forms
of $\chi^2_\p$)
can be automatically derived for non-starlike shape models or profiles
as well\cite{ipi}.

\section{Maximum compatibility estimate}

Let us choose as goodness-of-fit measures (from which probability
distributions can be constructed) the $\chi^2$-functions of $n$ data modes. 
Our task is to
construct a joint $\chi_{\rm tot}$ with well-defined weighting
for each data mode:
\be
\chi_{\rm tot}^2(P,D)=\chi_1^2(P,D_1)+\sum_{i=2}^n \lambda_{i-1}
\chi_i^2(P,D_i)\quad D=\{D_i,i=1,\ldots,n\}
\ee
(to which regularization functions $g(P)$ can be added), where $D_i$ denotes
the data from the source $i$, and 
$P\in\R^p$ is the set of model parameter values.
We assume the $\chi^2_i$-space to be nondegenerate, i.e.,
\[
\arg\min \chi^2_i(P) \ne \arg\min \chi^2_j(P),\quad i\ne j.
\]

In two dimensions, denote 
\bea
x(\lambda)&:=&\lbrace\chi_1^2\vert\min\chi_{\rm tot}^2;\lambda
\rbrace,\\\nonumber
y(\lambda)&:=&\lbrace\chi_2^2\vert\min\chi_{\rm tot}^2;\lambda\rbrace. 
\eea
The curve
\be
{\cal S}(\lambda):=[\log x(\lambda),\log y(\lambda)]
\ee
resembles the well-known ``L-curve''
related to, e.g., Tikhonov regularization \cite{belge,kaipio}.
However, here we make no assumptions on the shape of $\cal S$.
The curve $\cal S$ is a part of the boundary $\p\cal R$ of the region 
${\cal R}\in\R^2$ formed by the mapping $\chi: \R^p\rightarrow\R^2$
from the parameter space $\pp$ into $\chi_i^2$-space:
\[
\chi=\lbrace\pp\rightarrow
(\log\chi_1^2,\log\chi_2^2)\rbrace,\quad {\cal R}=\chi({\cal P})
\]
where the set ${\cal P}$ includes all the 
possible values of model parameters (assuming that $\chi$ is continuous and
well-behaved such that a connected $\cal R$ and $\p\cal R$ exist). 
If the possible values of $\chi^2_i$ are not bounded, the remaining part 
$\p\cal R\setminus\cal S$ stretches droplet-like towards $(\infty,\infty)$.
The parameter $\lambda$ describes the
position on the interesting part ${\cal S}\subset\p\cal R$,
and it is up to us to define a criterion
for choosing the optimal value of $\lambda$.
 
The logarithm
ensures that the shape of ${\cal S}(\lambda)$ is invariant under unit or scale
transforms in the $\chi_i^2$ as they merely translate ${\cal S}$ in the
$(\log\chi_1^2,\log\chi_2^2)$-plane. It also provides a meaningful metric
for the $\log\chi_i^2$-space: distances depict the relative difference
in $\chi^2$-sense, removing the problem of comparing the absolute values of 
quite different types of $\chi_i^2$.
The endpoints of ${\cal S}(\lambda)$ are at
$\lambda=0$ and $\lambda=\infty$, i.e., at the values of $\chi_i^2$
that result from using only one of the data modes in inversion. We can
translate the origin of the 
$(\log\chi_1^2,\log\chi_2^2)$-plane to a more natural
position by choosing the new coordinate axes to pass through these endpoints.
Denote
\bea
\hat x_0&=&\log x(\lambda)\vert_{\lambda=0}=\log\min\chi_1^2\\\nonumber
\hat y_0&=&\log y(\lambda)\vert_{\lambda\rightarrow\infty}=\log\min\chi_2^2. 
\eea
Then the ``ideal point'' $(\hat x_0,\hat y_0)$ is the new origin in the
$(\log x,\log y)$-plane. A natural choice for an optimal location
on $\cal S$ is the point closest to $(\hat x_0,\hat y_0)$, i.e., the
parameter values $P_0\in\pp$ such that
\be
P_0=\arg\min\Big([\log\chi^2_1(P)-\hat x_0]^2+[\log\chi^2_2(P)-
\hat y_0]^2\Big),
\ee
so we have, with $\lambda$ as argument,
\be
\lambda_0=\arg\min \Big([\log x(\lambda)-\hat x_0]^2+[\log y(\lambda)-
\hat y_0]^2\Big).\label{lambda0}
\ee
In this approach, neither the numbers of data points
in each $\chi^2_i$ nor the noise levels as such affect the
solution for the optimal $P_0$ as their scaling effects cancel out in
each quadratic term. $P_0$ is thus a pure compatibility estimate describing
the best model compromise explaining the datasets of different modes
simultaneously.

We call the point $P_0$ the
{\em maximum compatibility estimate} (MCE), and $\lambda_0$ the {\em maximum
compatibility weight} (MCW). This corresponds to the maximum
likelihood estimate in the case of one data mode, or to the maximum a
posteriori estimate as well since we can
include regularization functions here. If regularizing is used, the weights
for the functions are either determined in a similar manner (see below),
or they can be fixed and the regularization terms are absorbed in
$\chi_1^2$ (otherwise ${\cal S}\subset\p\cal R$ does not hold).

Another choice, frequently used in the L-curve approach, is to find the
$\lambda$ at which $\cal S$ attains its maximum curvature \cite{belge,kaipio}, 
but evaluating this point is less robust than finding $\lambda_0$,
and (\ref{lambda0}) is a more natural prescription, requiring no assumptions
on the shape of $\cal S$. We make two implicit assumptions here:
\begin{enumerate}
\item
The solutions $P_{\p\cal R}$ 
corresponding to points on $\p\cal R$ should be continuous (and one-to-one) 
in $\pp$-space along $\p\cal R$ at least in the vicinity of the solution
corresponding to $\lambda_0$. If this is not true (in practice, 
if $P_\lambda=\arg\min\chi^2_{\rm tot}(P)$
makes large jumps in $\pp$ for various $\lambda$ around $\lambda_0$), 
one should be cautious
about the uniqueness and stability of the chosen solution $P_0$, and
restrict the regions of $\pp$ included in the analysis. 
\item
The optimal point $\lambda_0$ on 
$\cal S$ should 
be feasible: if we have upper limits $\epsilon_i$ to acceptable $\chi^2_i$,
the feasible region $\cal F$
is the rectangle $\bigcap_i\lbrace\log\chi^2_i\le\log\epsilon_i\rbrace$. If
$[\log\chi_1^2(P_0),\log\chi_2^2(P_0)]\notin
\cal F$ and ${\cal F}\cap{\cal R}\ne
\emptyset$, we choose the point on the portion ${\cal S}\subset\cal R$ 
closest to the one corresponding
to $\lambda_0$ (i.e., $\log\chi^2_i=\log\epsilon_i$ for one $i$). 
If ${\cal F}\cap{\cal R}=\emptyset$, the data modes do not
allow a compatible joint model, so either the model is incorrect for 
one or both data modes, or one or both $\epsilon_i$ have been estimated 
too low (e.g., systematic errors have not been taken into account). Note
that model insufficiency should be taken into account in the estimation of 
$\epsilon_i$.
\end{enumerate}

Note that, in the interpretation ${\cal R}=\chi({\cal P})$,
$\lambda$, $\chi_{\rm tot}^2$ and $\p\cal R$ are all
in fact superfluous quantities, and we can locate the point estimate
MCE $P_0$ entirely without them with standard optimization procedures
(and with no extra computational cost).
However, it is useful (though computationally
somewhat noisier) to approximate $\cal S$ via the minimization 
of $\chi_{\rm tot}^2$ with sample values of $\lambda$ (see Fig.\ 1a), 
as in addition to obtaining the MCW $\lambda_0$ (and hence MCE as well) we can
plot $\cal S$ to examine the mutual behaviour of the complementary data
sources (including the position of the feasibility region $\cal F$ 
w.r.t. $\cal S$). The solution for $\lambda_0$ is also needed for
constructing distributions based on $\chi^2_{\rm tot}$. Another possibility
to examine $\cal R$ and $\p\cal R$ is direct adaptive Monte Carlo sampling, 
but this is computationally slow.

This approach straightforwardly generalizes to $n$ $\chi^2$-functions and
$n-1$ parameters $\lambda_i$ describing the position on the $n-1$-dimensional
boundary surface $\p\cal R$ of an $n$-dimensional domain $\cal R$: the MCE is
\be
P_0=\arg\min\sum_{i=1}^n\Big[\log\frac{\chi^2_i(P)}{\chi^2_{i0}}\Big]^2,
\quad\chi^2_{i0}:=\min\chi^2_i(P),\label{p0}
\ee
and the MCW is, for $\lambda\in\R^{n-1}$,
\be
\lambda_0=\arg\min\sum_{i=1}^n\Big[\log
\frac{\hat\chi_{i,{\rm tot}}^2(\lambda)}{\chi^2_{i0}}\Big]^2,
\quad\hat\chi_{i,{\rm tot}}^2(\lambda):=
\Big\{\chi^2_i\Big\vert\min\chi^2_{\rm tot};\lambda\Big\}.\label{l0}
\ee

Another scale-invariant version of MCE can be constructed by plotting 
$\chi^2_i$ in units of $\chi^2_i/\chi^2_{i0}$ and shifting the new origin
to $\chi^2_i/\chi^2_{i0}=1$:
\be
P_0=\arg\min\sum_{i=1}^n\Big[\frac{\chi^2_i(P)}{\chi^2_{i0}}-1\Big]^2,
\quad\lambda_0=\arg\min\sum_{i=1}^n\Big[
\frac{\hat\chi_{i,{\rm tot}}^2(\lambda)}{\chi^2_{i0}}-1\Big]^2.
\ee

This, however, is exactly the first-order approximation of (\ref{p0}) and
(\ref{l0}) in $\delta\ll 1$ when $\chi^2_i/\chi^2_{i0}=1+\delta$, giving
virtually the same result as (\ref{p0}) and (\ref{l0}) as usually
$\chi^2_i(P_0)/\chi^2_{i0}-1\ll 1$ in the region around $\chi^2_i(P_0)$, 
and any larger ratios of $\chi^2_i/\chi^2_{i0}$ are not eligible for
the optimal solution (see Fig.\ 1a).

It is possible to use this approach for general regularizing functions 
$g(P)$ as well (change $\chi_i^2\rightarrow g(P)$ for some $i$), 
but in such cases the shape of $\cal S$ must be taken into account.
If it is possible to have a solution $g(P')=0$ for a regularizing function
$g$ (or an almost vanishing $g(P')$ such that 
$\log g(P')\rightarrow -\infty$), 
one should, e.g., set a lower practical limit to 
$g(P)$ by looking at the shape of $S$, and choose the $\lambda_0$ within
the restricted part of $S$.

\section{Numerical implementation}

As examples of the optimal combining of brightness values and  
profile contours, we show
some results for asteroid data. 

Fig.\ 1a depicts a typical evaluation
of the curve $\cal S$ for 2 Pallas
at various choices of $\lambda$ (or rather, this plot
portrays the cross-section of the 2-surface $\p\cal R$ in $\R^3$ with
smoothness regularization weight fixed at its final optimal value).
The values for $\chi_i^2$ are normalized to be the rms deviations
of model fits $d_i=\sqrt{\chi_i^2/N_i}$, as in logarithmic scale
this corresponds only to a shift of origin and a uniform linear change
of plot scaling. The plotted points outline 
the curve ${\cal S}(\lambda)$ 
that is rather an oblique line than an L-shape, and the ideal
point region, i.e., the point closest to the lower left-hand corner,
can directly be found. The endpoints $\lambda=0$ and
$\lambda=\infty$ stop at saturation regions rather than continue to
large distances in the $\log\chi^2$-space.
As can be seen from Fig.\ 1a, 
computational noise in the estimated points at
$\lambda=0$ and $\lambda=\infty$, corresponding to a small change of the 
position of the new origin w.r.t. $\cal S$, 
does not affect the estimated location
of the optimal point on $\cal S$ significantly.

\begin{figure}
\psfig{file=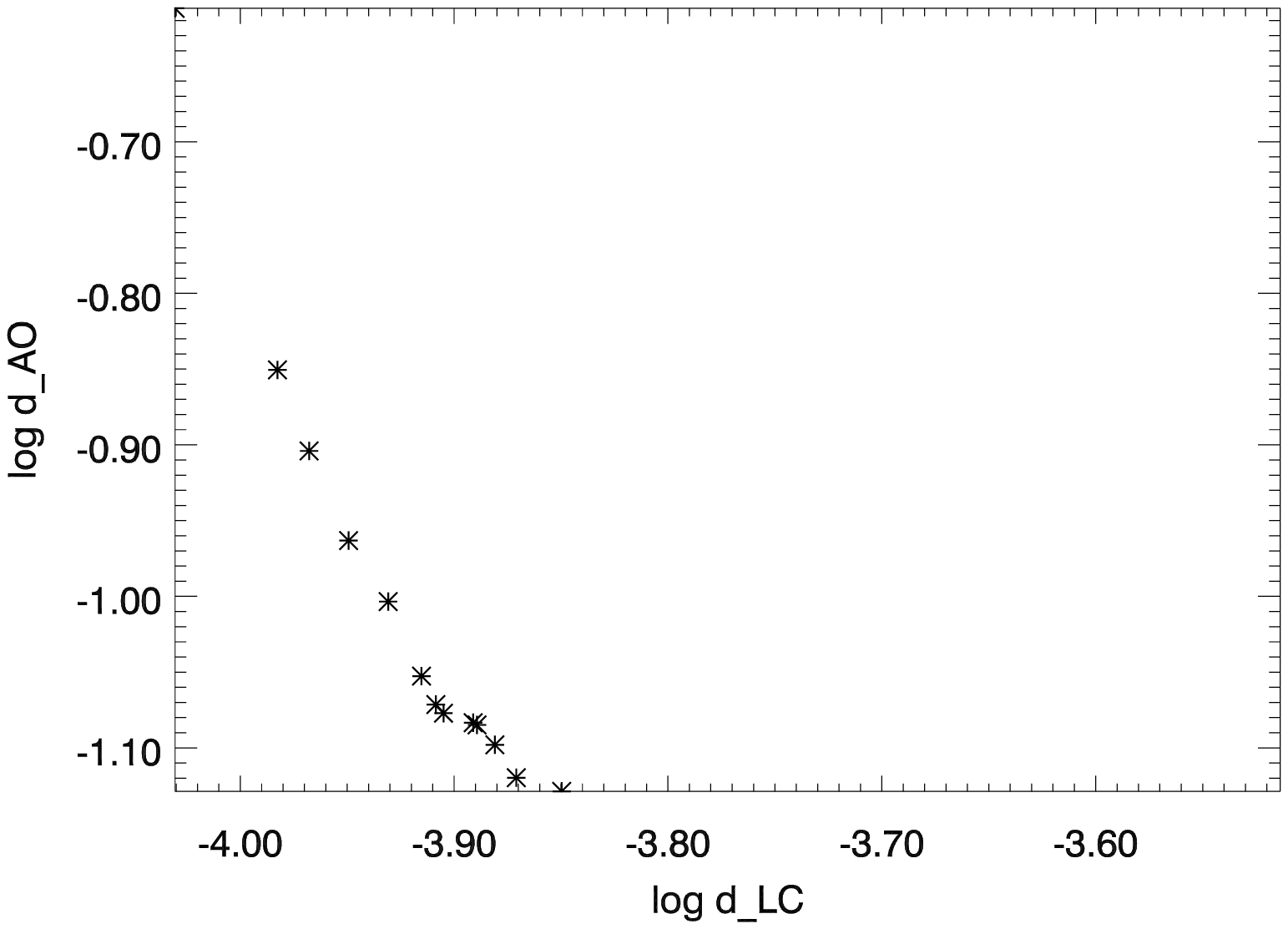,width=6cm}\psfig{file=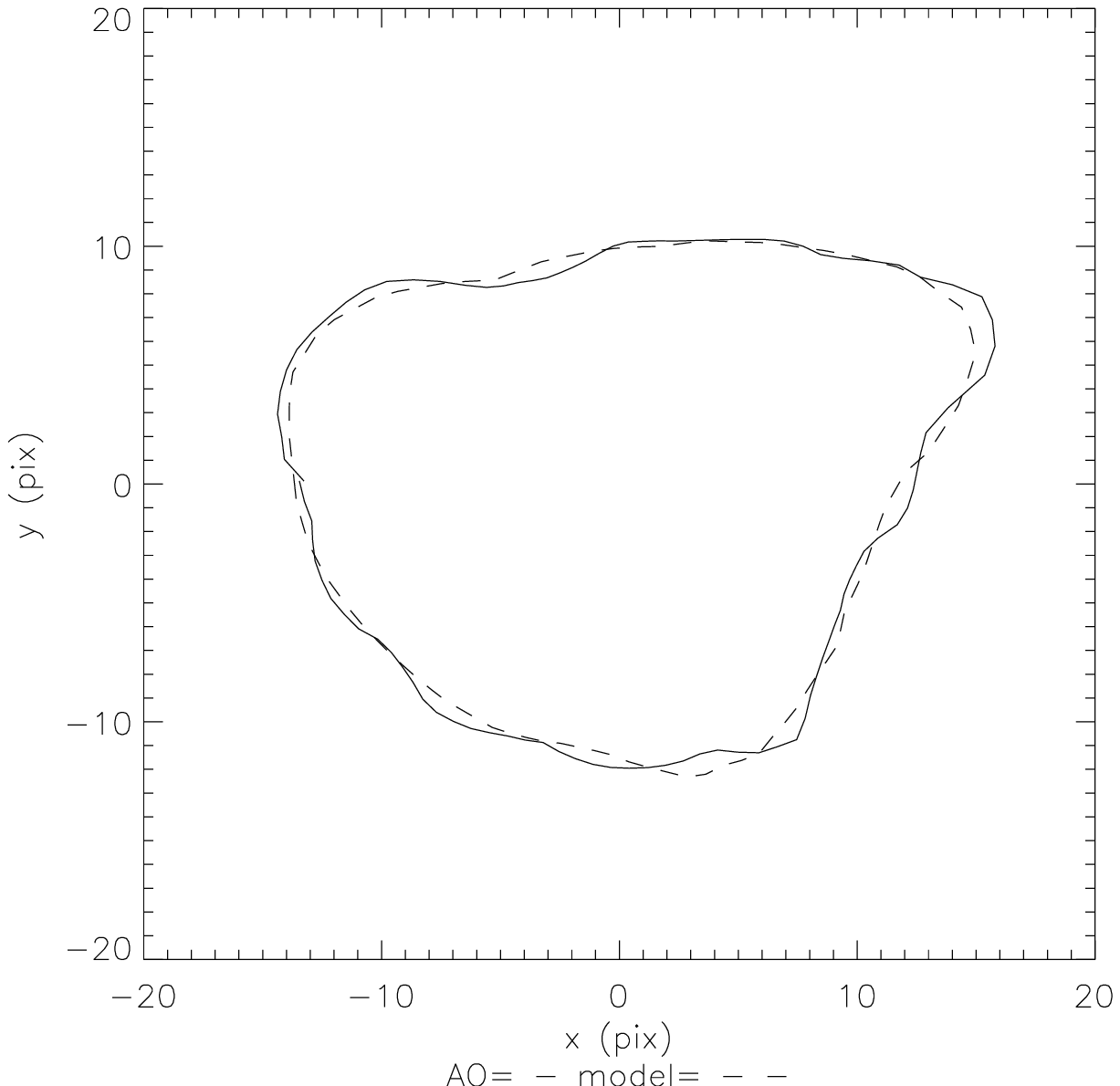,width=6cm}
\caption{(a) $\cal S$ curve plotted for 2 Pallas with various weights $\lambda$
(LC for brightness data, AO for adaptive optics profiles).
(b) Sample observed (solid line) vs. modelled (dashed line)
AO profile contour for 41 Daphne. Coordinates are in pixel units.}
\label{scurve}
\end{figure}

A sample observed vs.\ modelled profiles for 41 Daphne is shown
in Fig.\ 1b.
The starlike surface model was described by the exponential
Laplace (spherical harmonics) series for the surface radius $r$ (for
explicit positivity): 
\be
r(\theta,\varphi)=\exp\Big[\sum_{lm} c_{lm}Y_l^m(\theta,\varphi)\Big],
\quad (\theta,\varphi)\in S^2,
\ee
truncated at $l=8,m=6$, with $c_{lm}$ as the shape parameters 
to be solved for. Other model parameters are the profile offset 
$\varkappa_0$ for each image and the physical spin parameters
describing the rotational transformations yielding the correct
viewing and illumination directions $(\omega,\omega_0)$ on the body
\cite{genproj,ipi,yorp}.

\section{Discussion}

The concept of the maximum compatibility estimate is directly applicable
to any inverse problems with complementary data modes. The invariance 
properties of the MCE make it more generally usable than heuristic strategies
for choosing the weights, especially when they use assumptions on the
shape of $\p\cal R$ or other case-specific characteristics.
In our case study, the use of profiles is practical as it removes two 
sources of systematic
errors inherent to using full images (brightness distributions $\cal I$ on the
image plane): the errors in $\cal I$ from adaptive optics 
deconvolution and the model
$\cal I$ errors due to the insufficently modellable light-scattering 
properties of the surface of the target body.

\end{document}